\documentclass[12pt]{amsart}
\usepackage{amssymb}

\newcommand{\Bgp}{{\Z^\N}}

\long\def\forget#1\forgotten{}
\newcommand{\issuenumber}{25}
\newcommand{\issuemonth}{August}
\newcommand{\issueyear}{2008}

\newtheorem{issue}{Issue}

\theoremstyle{definition}

\theoremstyle{remark}

\newcommand{\ed}{
\newpage

\section{Unsolved problems from earlier issues}

\begin{issue}
Is $\binom{\Omega}{\Gamma}=\binom{\Omega}{\Tau}$?
\end{issue}

\begin{issue}
Is $\ufin(\cO,\Omega)=\sfin(\Gamma,\Omega)$?
And if not, does $\ufin(\cO,\Gamma)$ imply
$\sfin(\Gamma,\Omega)$?
\end{issue}

\stepcounter{issue}

\begin{issue}
Does $\sone(\Omega,\Tau)$ imply $\ufin(\Gamma,\Gamma)$?
\end{issue}

\begin{issue}
Is $\fp=\fp^*$? (See the definition of $\fp^*$ in that issue.)
\end{issue}

\begin{issue}
Does there exist (in ZFC) an uncountable set satisfying $\sfin(\B,\B)$?
\end{issue}

\stepcounter{issue}

\begin{issue}
Does $X \nin \NON(\M)$ and $Y\nin\mathsf{D}$ imply that
$X\cup Y\nin \COF(\M)$?
\end{issue}

\begin{issue}[CH]
Is $\split(\Lambda,\Lambda)$ preserved under finite unions?
\end{issue}

\begin{issue}
Is $\cov(\M)=\fo$? (See the definition of $\fo$ in that issue.)
\end{issue}

\begin{issue}
Does $\sone(\Gamma,\Gamma)$ always contain an element of cardinality $\fb$?
\end{issue}

\begin{issue}
Could there be a Baire metric space $M$ of weight $\aleph_1$ and a partition
$\mathcal{U}$ of $M$ into $\aleph_1$ meager sets where for each ${\mathcal U}'\subset\mathcal U$,
$\bigcup {\mathcal U}'$ has the Baire property in $M$?
\end{issue}

\stepcounter{issue} 

\begin{issue}
Does there exist (in ZFC) a set of reals $X$ of cardinality $\fd$ such that all
finite powers of $X$ have Menger's property $\sfin(\cO,\cO)$?
\end{issue}

\begin{issue}
Can a Borel non-$\sigma$-compact group be generated by a Hurewicz subspace?
\end{issue}

\begin{issue}[MA]
Is there an uncountable $X\sbst\R$ satisfying $\sone(\BO,\BG)$?
\end{issue}

\begin{issue}[CH]
Is there a totally imperfect $X$ satisfying $\ufin(\cO,\Gamma)$
that can be mapped continuously onto $\Cantor$?
\end{issue}

\begin{issue}[CH]
Is there a Hurewicz $X$ such that $X^2$ is Menger but not Hurewicz?
\end{issue}

\begin{issue}
Does the Pytkeev property of $C_p(X)$ imply that $X$ has Menger's property?
\end{issue}

\begin{issue}
Does every hereditarily Hurewicz space satisfy $\sone(\BG,\BG)$?
\end{issue}

\begin{issue}[CH]
Is there a Rothberger-bounded $G\le\Bgp$ such that $G^2$ is not Menger-bounded?
\end{issue}

\begin{issue}
Let $\cW$ be the van der Waerden ideal.
Are $\cW$-ultrafilters closed under products?
\end{issue}

\begin{issue}
Is the $\delta$-property equivalent to the $\gamma$-property $\binom{\Omega}{\Gamma}$?
\end{issue}

\stepcounter{issue}

\stepcounter{issue}

\general\end{document}}

\newcommand{\Cantor}{{\{0,1\}^\N}}

\newcommand{\fb}{\mathfrak{b}}

\newcommand{\fc}{\mathfrak{c}}
\newcommand{\fd}{\mathfrak{d}}
\newcommand{\fp}{\mathfrak{p}}

\newcommand{\NON}{{\mathsf   {NON}}}
\newcommand{\COF}{{\mathsf   {COF}}}

\newcommand{\M}{\mathcal{M}}

\newcommand{\op}{\operatorname}
\newcommand{\cov}{\mathsf{cov}}

\newcommand{\R}{\mathbb{R}}
\newcommand{\Q}{\mathbb{Q}}

\newcommand{\fo}{\mathfrak{od}}

\newcommand{\w}{\omega}

\renewcommand{\split}{\mathsf{Split}}
\newcommand{\bq}{\begin{quote}}
\newcommand{\eq}{\end{quote}}
\newcommand{\cO}{\mathcal{O}}
\newcommand{\B}{\mathcal{B}}
\newcommand{\BG}{\B_\Gamma}

\newcommand{\BO}{\B_\Omega}

\newcommand{\sone}{\mathsf{S}_1}    \newcommand{\sfin}{\mathsf{S}_\mathrm{fin}}

\newcommand{\ufin}{\mathsf{U}_\mathrm{fin}}

\newcommand{\nin}{\not\in}

\newcommand{\cW}{\mathcal{W}}

\newcommand{\N}{\mathbb{N}}
\newcommand{\Z}{\mathbb{Z}}

\newcommand{\sbst}{\subseteq}
\newcommand{\by}[2]{\par\hfill\emph{#1}, #2}
\newcommand{\nby}[1]{\par\hfill\emph{#1}}
\newcommand{\Tau}{\mathrm{T}}
\newcommand{\CE}{\textsc{CE}}

\newcommand{\be}{\begin{enumerate}}
\newcommand{\ee}{\end{enumerate}}
\newcommand{\bi}{\begin{itemize}}
\newcommand{\ei}{\end{itemize}}

\newcommand{\general}{\small\vfill\par\noindent\hrulefill\par
\noindent\textbf{Previous issues.} The previous issues of this
bulletin are available online at\\
\texttt{http://front.math.ucdavis.edu/search?\&t=\%22SPM+Bulletin\%22}
\\[0.1cm]
\textbf{Contributions.} Announcements, discussions, and open problems should be emailed
to \texttt{tsaban@math.biu.ac.il}\\[0.1cm]
\textbf{Subscription.}
To receive this bulletin (free) to your e-mailbox, e-mail us.
}

\newcommand{\link}[1]{\par\hfill{\texttt{#1}}}

\newcommand{\arXiv}[5]{\subsection{#2}{#4}\par\hfill{\arx{#1}}\par\hfill\emph{#3}}

\newcommand{\nAMSPaper}[4]{\subsection{#2}{#4}\par\hfill{\texttt{#1}}\par\hfill\emph{#3}}

\newcommand{\arx}[1]{\texttt{http://arxiv.org/abs/#1}}
\newcommand{\url}[1]{\bq\texttt{#1}\eq}
\newcommand{\online}[1]{The paper is available online at \url{#1}}

\title[$\mathcal{SPM}$ Bulletin \textbf{\issuenumber} (\issuemonth{} \issueyear)]{%
$\mathcal{SPM}$ Bulletin\\[0.5cm]
Issue number \issuenumber: \issuemonth{} \issueyear{} \CE{}}

\begin{document}
\maketitle

\tableofcontents

\section{Editor's note}

The slides of the talks given at the conference \emph{Ultramath} (Pisa, Italy, June 2008)
are available at
\url{http://www.dm.unipi.it/\~{}ultramath/abstracts.html}

Enjoy.

\medskip

\by{Boaz Tsaban}{tsaban@math.biu.ac.il}

\hfill \texttt{http://www.cs.biu.ac.il/\~{}tsaban}

\section{Research announcements}

\arXiv{0803.3498}
{Combinatorial and model-theoretical principles related to regularity of ultrafilters and compactness of topological spaces, I}
{Paolo Lipparini}
{We begin the study of the consequences of the existence of certain infinite
matrices. Our present application is to compactness of products of topological
spaces.}

\arXiv{0803.4117}
{Fr\'{e}chet-Urysohn fans in free topological groups}
{Taras Banakh, Du\v{s}an Repov\v{s}, and Lyubomyr Zdomskyy}
{In this paper we answer the question of T. Banakh and M. Zarichnyi
constructing a copy of the Fr\'echet-Urysohn fan $S_\omega$ in a topological group
$G$ admitting a functorial embedding $[0,1]\subset G$. The latter means that
each autohomeomorphism of $[0,1]$ extends to a continuous homomorphism of $G$.
This implies that many natural free topological group constructions (e.g. the
constructions of the Markov free topological group, free abelian topological
group, free totally bounded group, free compact group) applied to a Tychonov
space $X$ containing a topological copy of the space $\Q$ of rationals give
topological groups containing $S_\omega$.}

\arXiv{0804.1333}
{Packing index of subsets in Polish groups}
{Taras Banakh, Nadya Lyaskovska, and Du\v{s}an Repov\v{s}}
{For a subset $A$ of a Polish group $G$, we study the (almost) packing index
$\op{ind}_P(A)$ (resp. $\op{Ind}_P(A)$) of $A$, equal to the supremum of cardinalities
$|S|$ of subsets $S\subset G$ such that the family of shifts $\{xA\}_{x\in S}$
is (almost) disjoint (in the sense that $|xA\cap yA|<|A|$ for any distinct
points $x,y\in S$). Subsets $A\subset G$ with small (almost) packing index are
small in a geometric sense. We show that $\op{ind}_P(A)\in \N\cup\{\aleph_0,\fc\}$
for any $\sigma$-compact subset $A$ of a Polish group. If $A\subset G$ is
Borel, then the packing indices $\op{ind}_P(A)$ and $\op{Ind}_P(A)$ cannot take values
in the half-interval $[\mathfrak{sq}(\Pi^1_1),\fc)$ where $\mathfrak{sq}(\Pi^1_1)$ is a certain
uncountable cardinal that is smaller than $\fc$ in some models of ZFC. In each
non-discrete Polish Abelian group $G$ we construct two closed subsets
$A,B\subset G$ with $\op{ind}_P(A)=\op{ind}_P(B)=\fc$ and $\op{Ind}_P(A\cup B)=1$ and then
apply this result to show that $G$ contains a nowhere dense Haar null subset
$C\subset G$ with $\op{ind}_P(C)=\op{Ind}_P(C)=\kappa$ for any given cardinal number
$\kappa\in[4,\fc]$.}

\arXiv{0804.1335}
{Symmetric monochromatic subsets in colorings of the Lobachevsky plane}
{T. Banakh, A. Dudko and D. Repov\v{s}}
{We prove that for each partition of the Lobachevsky plane into finitely many
Borel pieces one of the cells of the partition contains an unbounded centrally
symmetric subset.}

\arXiv{0804.1593}
{Structural Ramsey theory of metric spaces and topological dynamics of isometry groups}
{L. Nguyen Van Th\'e}
{In 2003, Kechris, Pestov and Todorcevic showed that the structure of certain
separable metric spaces - called ultrahomogeneous - is closely related to the
combinatorial behavior of the class of their finite metric spaces. The purpose
of the present paper is to explore the different aspects of this connection.}

\arXiv{0804.4019}
{Distinguishing Number of Countable Homogeneous Relational Structures}
{C. Laflamme, L. Nguyen Van Th\'e, N. W. Sauer}
{The distinguishing number of a graph $G$ is the smallest positive integer $r$
such that $G$ has a labeling of its vertices with $r$ labels for which there is
no non-trivial automorphism of $G$ preserving these labels. Albertson and
Collins computed the distinguishing number for various finite graphs, and
Imrich, Klav\v{z}ar and Trofimov computed the distinguishing number of some
infinite graphs, showing in particular that the Random Graph has distinguishing
number 2. We compute the distinguishing number of various other finite and
countable homogeneous structures, including undirected and directed graphs, and
posets. We show that this number is in most cases two or infinite, and besides
a few exceptions conjecture that this is so for all primitive homogeneous
countable structures.}

\arXiv{0804.4548}
{Indestructible colourings and rainbow Ramsey theorems}
{Lajos Soukup}
{We give a negative answer to a question of Erdos and Hajnal: it is consistent
that GCH holds and there is a colouring $c:[{\omega_2}]^2\to 2$
establishing $\omega_2 \not\to [(\omega_1;{\omega})]^2_2$ such
that some colouring $g:[\omega_1]^2\to 2$ can not be embedded into
$c$. It is also consistent that $2^{\omega_1}$ is arbitrarily
large, and a function $g$ establishes $2^{\omega_1} \not\to
[(\omega_1,\omega_2)]^2_{\omega_1}$ such that there is no
uncountable $g$-rainbow subset of $2^{\omega_1}$. We also show
that for each $k\in {\omega}$ it is consistent with Martin's Axiom
that the negative partition relation $\omega_1 \not\to^*
[(\omega_1;\omega_1)]_{k-bdd}$ holds.}

\nAMSPaper{www.ams.org/proc/0000-000-00/S0002-9939-08-09334-9}
{Products of Borel subgroups}
{Longyun Ding and Bingqing Li}
{We investigate the Borelness of the product of two Borel
subgroups in Polish groups. While the intersection of these two
subgroups is Polishable, the Borelness of their product is
confirmed. On the other hand, we construct two $\Delta^0_3$ subgroups whose
product is not Borel in every uncountable abelian Polish group. }

\nAMSPaper{http://www.ams.org/journal-getitem?pii=S0002-9939-08-09548-8}
{Selection theorems and treeability}
{Greg Hjorth}
{We show that domains of non-trivial $\Sigma^1_1$ trees have $\Delta^1_1$ members.
Using this, we show that smooth treeable equivalence relations
have Borel transversals, and essentially countable treeable
equivalence relations have Borel complete countable sections. We
show also that treeable equivalence relations which are ccc
idealistic, measured, or generated by a Borel action of a Polish
group have Borel complete countable sections.}

\arXiv{0805.1548}
{Combinatorial and model-theoretical principles related to regularity of ultrafilters and compactness of topological spaces, IV}
{Paolo Lipparini}
{We extend to singular cardinals the model-theoretical relation
$\lambda\stackrel{\kappa}{\Rightarrow} \mu$ introduced in P. Lipparini, The compactness
spectrum of abstract logics, large cardinals and combinatorial principles,
Boll. Unione Matematica Italiana ser. VII, {\bf 4-B} 875--903 (1990). We extend
some results obtained in Part II, finding equivalent conditions involving
uniformity of ultrafilters and the existence of certain infinite matrices. Our
present definition suggests a new compactness property for abstract logics.}

\arXiv{0806.2719}
{A property of $C_p[0,1]$}
{Michael Levin}
{We prove that for every finite dimensional compact metric space $X$ there is
an open continuous linear surjection from $C_p[0,1]$ onto $C_p(X)$. The proof
makes use of embeddings introduced by Kolmogorov and Sternfeld in connection
with Hilbert's 13th problem.}

\subsection{A Dedekind Finite Borel Set}
In this paper we prove three theorems about the theory of Borel sets in
models of ZF without any form of the axiom of choice. We prove that if $B$ is a
$G_{\delta\sigma}$ set, then either $B$ is countable or $B$ contains a perfect subset.
Second, we prove that if the real line is the countable union of countable
sets, then there exists an $F_{\sigma\delta}$ set which is uncountable but contains
no perfect subset. Finally, we construct a model of ZF in which we have an
infinite Dedekind finite set of reals which is $F_{\sigma\delta}$.
\par\hfill{\texttt{http://www.math.wisc.edu/\~{}miller/res/ded.pdf}}
\par\hfill\emph{Arnold W. Miller}

\arXiv{0806.4499}
{Aronszajn Compacta}
{Joan E. Hart and Kenneth Kunen}
{We consider a class of compacta $X$ such that the maps from $X$ onto metric
compacta define an Aronszajn tree of closed subsets of $X$.}

\arXiv{0806.4220}
{A strong antidiamond principle compatible with CH}
{James Hirschorn}
{A strong antidiamond principle $(*c)$ is shown to be consistent with CH. This
principle can be stated as a ``$P$-ideal dichotomy'': every $P$-ideal on $omega-1$
(i.e. an ideal that is $\sigma$-directed under inclusion modulo finite) either has
a closed unbounded subset of $\omega_1$ locally inside of it, or else has a
stationary subset of $\omega_1$ orthogonal to it. We rely on Shelah's theory of
parameterized properness for NNR iterations, and make a contribution to the
theory with a method of constructing the properness parameter simultaneously
with the iteration. Our handling of the application of the NNR iteration theory
involves definability of forcing notions in third order arithmetic, analogous
to Souslin forcing in second order arithmetic.}

\arXiv{0806.4732}
{On the strength of Hausdorff's gap condition}
{James Hirschorn}
{Hausdorff's gap condition was satisfied by his original 1936 construction of
an $(\omega_1,\omega_1)$ gap in $P(\N)/Fin$. We solve an open problem in determining
whether Hausdorff's condition is actually stronger than the more modern
indestructibility condition, by constructing an indestructible
$(\omega_1,\omega_1)$ gap not equivalent to any gap satisfying Hausdorff's
condition, from uncountably many random reals.}

\arXiv{0807.0147}
{Nonhomogeneous analytic families of trees}
{James Hirschorn}
{We consider a dichotomy for analytic families of trees stating that either
there is a colouring of the nodes for which all but finitely many levels of
every tree are nonhomogeneous, or else the family contains an uncountable
antichain. This dichotomy implies that every nontrivial Souslin poset
satisfying the countable chain condition adds a splitting real.
 We then reduce the dichotomy to a conjecture of Sperner Theory. This
conjecture is concerning the asymptotic behaviour of the product of the sizes
of the m-shades of pairs of cross-t-intersecting families.}

\arXiv{0806.4760}
{Reasonable non-Radon-Nikodym ideals}
{Vladimir Kanovei and Vassily Lyubetsky}
{For any abelian Polish $\sigma$-compact group $H$ there exist a $\sigma$-ideal $Z$ over
$\N$ and a Borel $Z$-approximate homomorphism $f : H \to H^\N$ which is not
$Z$-approximable by a continuous true homomorphism $g : H \to H^\N$.}

\arXiv{0807.1254}
{$\sigma$-continuity and related forcings}
{Marcin Sabok}
{The Steprans forcing notion arises as a quotient of Borel sets modulo the
ideal of $\sigma$-continuity of a certain Borel not $\sigma$-continuous
function. We give a characterization of this forcing in the language of trees
and using this characterization we establish such properties of the forcing as
fusion and continuous reading of names. Although the latter property is usually
implied by the fact that the associated ideal is generated by closed sets, we
show it is not the case with Steprans forcing. We also establish a connection
between Steprans forcing and Miller forcing thus giving a new description of
the latter. Eventually, we exhibit a variety of forcing notions which do not
have continuous reading of names in any presentation.}

\arXiv{0807.2205}
{An exact Ramsey principle for block sequences}
{Christian Rosendal}
{We prove an exact, i.e., formulated without $\Delta$-expansions, Ramsey
principle for infinite block sequences in vector spaces over countable fields,
where the two sides of the dichotomic principle are represented by respectively
winning strategies in Gowers' block sequence game and winning strategies in the
infinite asymptotic game. This allows us to recover Gowers' dichotomy theorem
for block sequences in normed vector spaces by a simple application of the
basic determinacy theorem for infinite asymptotic games.}

\nAMSPaper{http://www.ams.org/tran/0000-000-00/S0002-9947-08-04503-0/home.html}
{Baire reflection}
{Stevo Todorcevic and Stuart Zoble}
{We study reflection principles involving nonmeager sets and the
Baire Property which are consequences of the generic
supercompactness of $\w_2$, such as the principle asserting that any
point countable Baire space has a stationary set of closed
subspaces of weight $\w_1$ which are also Baire spaces. These principles
entail the analogous principles of stationary reflection but are
incompatible with forcing axioms. Assuming MM, there is a Baire
metric space in which a club of closed subspaces of weight $\w_1$ are
meager in themselves. Unlike stronger forms of Game Reflection,
these reflection principles do not decide CH, though they do give $\w_2$
as an upper bound for the size of the continuum.}

\arXiv{0807.3978}
{Tukey classes of ultrafilters on $\w$}
{David Milovich}
{Motivated by a question of Isbell, we show that Jensen's Diamond Principle
implies there is a non-P-point ultrafilter U on $\w$ such that U, whether
ordered by reverse inclusion or reverse inclusion mod finite, is not Tukey
equivalent to the finite sets of reals ordered by inclusion. We also show that,
for every regular infinite kappa not greater than $2^{\aleph_0}$, if
$MA(\sigma-centered)$ holds, then some ultrafilter U on $\w$, ordered by reverse
inclusion mod finite, is Tukey equivalent to the sets of reals of size less
than kappa, ordered by inclusion. We also prove two negative ZFC results about
the possible Tukey classes of ultrafilters on $\w$.}

\arXiv{0807.3846}
{Countably determined compact abelian groups}
{Dikran Dikranjan, Dmitri Shakhmatov}
{For an abelian topological group $G$ let $\widehat{G}$ be the dual
group of all continuous characters endowed with the compact open
topology.  Given a closed subset $X$ of an infinite compact abelian
group $G$ such that $w(X)<w(G)$ and an open neighbourhood $U$ of $0$
in $T$, we show that $|\{\pi\in\widehat{G}: \pi(X)\subseteq
U\}|=|\widehat{G}|$.  (Here $w(G)$ denotes the weight of $G$.)  A
subgroup $D$ of $G$ determines $G$ if the restriction homomorphism
$\widehat{G}\to \widehat{D}$ of the dual groups is a topological
isomorphism.  We prove that $w(G)=\min\{|D|: D$ is a subgroup of $G$
that determines $G\}$ for every compact abelian group $G$.  In
particular, an infinite compact abelian group determined by its
countable subgroup must be metrizable.  This gives a negative answer
to questions of Comfort, Hern\' andez, Macario, Raczkowski and
Trigos-Arrieta.  As an
application, we furnish a short elementary proof of the result
that compact determined abelian groups are metrizable.}

\subsection{A topological reflection principle equivalent to Shelah's Strong Hypothesis}
We notice that Shelah's Strong Hypothesis is equivalent to the following reflection principle:
Suppose $\langle X,\tau\rangle$ is a first-countable space whose density is a regular cardinal, $\kappa$.
If every  separable subspace of $X$ is of cardinality at most  $\kappa$, then the cardinality of $X$ is $\kappa$.
\link{dx.doi.org/10.1090/S0002-9939-08-09411-2}
\nby{Assaf Rinot}

\arXiv{0808.1654}
{Superfilters, Ramsey theory, and van der Waerden's Theorem}
{Nadav Samet and Boaz Tsaban}
{Superfilters are generalized ultrafilters, which capture the underlying
concept in Ramsey theoretic theorems such as van der Waerden's Theorem. We
establish several properties of superfilters, which generalize both Ramsey's
Theorem and its variant for ultrafilters on the natural numbers. We use them to
confirm a conjecture of Ko\v{c}inac and Di Maio, which is a generalization of a
Ramsey theoretic result of Scheepers, concerning selections from open covers.
Following Bergelson and Hindman's 1989 Theorem, we present a new simultaneous
generalization of the theorems of Ramsey, van der Waerden, Schur,
Folkman-Rado-Sanders, Rado, and others, where the colored sets can be much
smaller than the full set of natural numbers.}

\ed